\documentclass{amsart}

\usepackage[english]{babel} 
\usepackage[latin1]{inputenc}
\usepackage[T1]{fontenc} 
\usepackage{graphicx} 

\usepackage{verbatim} 

\usepackage{hyperref} 

\newcommand{\bbR}{\mathbb{R}}

\newcommand{\bbZ}{\mathbb{Z}}

\newcommand{\bbTP}{\mathbb{TP}}

\newcommand{\cL}{\mathcal{L}}

\newcommand{\val}{{\mathrm{val}}} 

\newcommand{\supp}{{\mathrm{Supp}}} 

\newcommand{\ord}{{\mathrm{ord}}} 

\newcommand{\Div}{{\mathrm{Div}}} 

\newcommand{\rk}{{\mathrm{rk}}} 

\newcommand{\f}{\varphi}

\newcommand{\G}{\Gamma}

\newcommand{\8}{\infty}
\newcommand{\0}{\emptyset}

\newcommand{\meno}{\smallsetminus}
\newcommand{\apa}{{``}}

\newcommand{\bs}{\bigskip}
\renewcommand{\ss}{\smallskip}

\usepackage{amssymb}
\usepackage{amsthm}

\newtheorem{Lemma}{Lemma}
\newtheorem{Proposition}{Proposition}

\newtheorem{Corollary}{Corollary}
\newtheorem{Theorem}{Theorem}

\begin{document}

\title{On tropical Clifford's Theorem}
\author{Laura Facchini}
\address{Department of Mathematics\\
University of Trento\\
Via Sommarive 14\\
38123 Povo (TN), Italy}
\email{laura.facchini@unitn.it}
\thanks{This work is part of the author's PhD research program, advised by Claudio Fontanari and partially supported by MIUR (Italy).}
\subjclass{14T05; 14H51}
\keywords{Tropical curve; linear system; Clifford's inequality; hyperelliptic curve}

\begin{abstract}
We prove an analogue of Clifford's inequality for tropical curves. Next we focus on the hyperelliptic case and we characterize divisors attaining equality. Finally we speculate whether inequality in tropical Clifford's Theorem does imply hyperellipticity, in analogy with the classical case.
\end{abstract}

\maketitle

\section{Introduction}

The geometry of tropical curves is a very hot topic, from the point of view both of Brill-Noether theory (see for instance \cite{CDPR}) and of moduli theory (see for instance \cite{Cap}).

In this short note, we investigate Clifford's Theorem for tropical curves. Indeed, we are aware of at least three independent proofs of the tropical Riemann-Roch Theorem (see \cite{HKN}, \cite{MiZh} and \cite{GaKe}), while to the best of our knowledge Clifford's inequality has been established up to now only for graphs by Baker and Norine (see \cite{BN1}, Corollary 3.5).

Here we extend Baker-Norine's result to tropical curves and we characterize divisors attaining equality in Clifford's Theorem at least for hyperelliptic curves (see Theorem \ref{clifford}).

For genus $g \leq 4$ these ones are the only cases (with the obvious exceptions of the zero and the canonical divisors), as in the classical setting. The standard proofs for curves $C$ of higher genus (see \cite{Har}, IV, Theorem 5.4, and \cite{ACGH}, p. 109) rely either on the linear algebra of the subspaces $\cL(D) := \{ f \in K(C) : (f) \geq -D \}$ or on the fact that if equality holds in the statement of Clifford's Theorem for a divisor $D \ne 0,K$, then the canonical divisor $K$ cannot be very ample and the curve must be hyperelliptic. Unluckily, in the tropical case $\cL(D)$ is definitely not a vector space (see for instance \cite{BN1}, Remark 1.13.(iii)) and there are examples of hyperelliptic curves with very ample canonical divisor (namely, a flower with at least three petals: see \cite{HMY} after the proof of Theorem 46 and notice that the equality in Clifford's Theorem is still attained by the multiples of the $g^1_2$ even if the canonical divisor is very ample). Hence there seems to be no hope to adapt the classical arguments to the tropical context.   

\section{Divisors on tropical curves}

First, we recall some definitions and remarks about divisors on tropical curves. 
For more details, see for instance \cite{MiZh}, \textsection 7. 

\ss

A \emph{multigraph} is a graph which is permitted to have multiple edges, i.e. edges that have the same end nodes.

With \emph{graph} we mean a finite and connected multigraph, not necessarily loop-free (i.e. there may be edges that connect a vertex to itself). 

A \emph{metric graph} is a pair $(\G,l)$ consisting of a graph $\G=(V(\G),E(\G))$ together with the so-called \emph{length function} $l : E(\G) \to \bbR_{> 0}$ which associates to every edge of the graph its length. 

Following \cite{GaKe}, Definition 1.1.(c), and \cite{HKN}, \textsection 1.1, by a \emph{tropical curve} we mean a metric graph with possibly unbounded ends, i.e. a pair $(\G,l)$ where the length function takes values in $\bbR_{> 0} \cup \{ \8 \}$. 

\ss

The tropical divisor group $\Div (\G)$ on a tropical curve $\G$ is defined as all formal finite linear combinations 
$$D=\sum a_i p_i$$ 
with $a_i \in \bbZ$ and $p_i \in \G$ for all $i$.

The \emph{support} of a divisor $D=\sum a_i p_i$ is defined as 
$$\supp D := \{ p_i \in \G : a_i \neq 0 \}.$$

We say that $D = \sum a_i p_i$ is \emph{effective} if $a_i \geq 0$ for all $i$ and we write $D \geq 0$.

For $D = \sum a_i p_i \in \Div(\G)$, we define the \emph{degree} of $D$ by the formula
$$\deg(D) = \sum a_i.$$

Following \cite{GaKe}, Definition 1.4, a \emph{rational function} on a tropical curve $\G$ is a continuous function $f : \G \to \bbR \cup \{\pm \8\}$ such that the restriction of $f$ to any edge of $\G$ is a piecewise linear integer affine function with a finite number of pieces. 

For a rational function $f$, the \emph{order} $\ord_P f \in \bbZ$ of a point $P \in \G$ is the sum of the slopes of all its outgoing segments.

As in \cite{GaKe}, we note that $\ord_P f = 0$ for all points $P \in \G \meno V(\G)$ at which $f$ is locally linear and thus for all but finitely many points. 

We can therefore define a \emph{principal divisor}, i.e. a divisor associated to a tropical rational function 
$$(f) := \sum_{P \in V(\G)} \ord_P f \cdot P \ \in \Div(\G).$$

As in the classical case, we can naturally define \emph{linearly equivalent divisors} as tropical divisors whose difference is a principal divisor.
As usual we denoted two linearly equivalent divisors with $D \sim E$.

\ss

For a divisor $D$, we use $|D|$ to denote the \emph{(complete) linear system} associated to it, i.e.
$$|D| = \{E \in \Div(\G) : E \geq 0, \ E \sim D \}.$$

\bs

We stress that the set
$$\cL(D) := \{ f \in K(\G) : (f) \geq -D \}$$
where $K(\G)$ is the function semifield on $\G$, is not a vector space, unlike the classical case (see \cite{BN1}, Remark 1.13.(iii)), but only a finitely generated semimodule (see \cite{MiZh}, Definition 2.1). Its minimal number of generators turns out to be a useful invariant, since a finite set $(f_0, \ldots, f_r)$ of generators for $\cL(D)$ induces a map
\begin{eqnarray*}
\f: \G & 	\to & 	\bbTP^{r-1} \\
	  x & \mapsto & (f_0(x):\ldots:f_r(x))
\end{eqnarray*}
(see \cite{HMY}, \textsection 5). As pointed out in \cite{HMY}, Lemma 38, an effective divisor $D$ is very ample (i.e., $\cL(D)$ separates points) if and only if $\f$ is injective for any choice of a set of generators for $\cL(D)$.

It turns out that $|D|$ is the tropical projectivization of $\cL(D)$ (see \cite{MiZh}, Definition 2.6 and \textsection 4.4) and it has the structure of a compact CW-complex (see \cite{MiZh}, \textsection 7). However, $|D|$ has not pure topological dimension (see \cite{GaKe}, Example 1.11) and its tropical dimension (which coincides with the topological dimension according to \cite{MiZh}, Proposition 2.5) does not satisfy the statement of Riemann-Roch (see \cite{BN1}, Remark 1.13 (iii)).
Thus, following Baker and Norine \cite{BN1}, we need to introduce the notion of rank.

\ss

Let $D$ be an effective divisor on $\G$. The \emph{rank} $\rk|D|$ of the linear system $|D|$ is defined as
\begin{equation} \label{rango}
\rk|D| := \max\{ r : |D-E| \neq \0 \quad \forall E \geq 0 \mbox{ with } \deg E = r \}
\end{equation}
and $\rk|D| := -1$ if $|D|=\0$.

In the classical case, $\rk|D|$ is called the dimension of the linear system $|D|$, but in our setting according to \cite{HKN}, \textsection 1.1, there is no known interpretation of $\rk|D|$ as the topological dimension of some space. Thus, we refer to $\rk|D|$ as \apa the rank'' rather than \apa the dimension'', following \cite{HKN}.

The definition of rank immediately implies the following fact (see for instance \cite{BN1}, Lemma 2.1): 

\begin{Lemma} \label{somma}
Let $D, E$ effective divisors on a tropical curve $\G$. Then
$$\rk|D|+\rk|E| \ \leq \ \rk|D+E|.$$
\end{Lemma}

\ss

We also recall from \cite{HMY}, \textsection 9, the following open question: is there any relation between $\rk|D|$ and the minimal number of generators of $\cL(D)$?

If we define the dimension of a semimodule as the maximum dimension of its straight reflexive submodules (see \cite{Yos}, \textsection 2), then we have
$$\rk|D| \leq \dim \cL(D)-1$$
(see \cite{Yos}, Theorem 2.7).

\bs

The \emph{canonical divisor} on $\G$ is the divisor $K$ defined as
$$K = \sum (\val(p_i) - 2) p_i,$$
where $\val(p_i)$ represents the valence of the point $p_i$, i.e. the number of its incident edges.

The \emph{genus} of $\G$ is given by Euler's formula, that states
$$g = |E(\G)| - |V(\G)| + 1,$$
where $|E(\G)|$ and $|V(\G)|$ represent respectively the number of edges and vertices of $\G$.

\ss

The tropical version of the Riemann-Roch Theorem (see \cite{HKN}, \cite{MiZh} and \cite{GaKe}) states that
$$\rk|D|-\rk|K-D| = \deg D -g+1.$$

As in the classical contest, we define a \emph{special divisor} as an effective divisor $D$ such that 
$$\rk|K-D| = \rk|D| - \deg D +g-1 \ \geq \ 0.$$

We say that a tropical curve is \emph{hyperelliptic} if there exists a divisor $D$ on it such that $|D|$ is a linear system of rank $1$ and degree $2$, i.e. a $g_2^1$.

In the tropical case it is still true that if the canonical divisor $K$ is not very ample then $\Gamma$ is hyperelliptic, 
but the converse turns out to be false (see \cite{HMY}, Theorem 46 and the example following it, and notice that this 
does not contradict \cite{BN2}, Proposition 5.22, where the so-called canonical map is not associated to the canonical 
divisor). 

As remarked by \cite{Bak}, Remark 3.4, in the case of metric graphs, we have:
 
\begin{Lemma} \label{g12}
If a $g_2^1$ exists, it is automatically unique.
\end{Lemma}

We also point out the following fact:

\begin{Proposition} \label{sistema}
Let $D$ an effective divisor on a tropical curve $\G$ of genus $g \geq 2$. Suppose that $\deg D =2g-2$ and $\rk|D| \geq g-1$. Then $|K|=|D|$.
\end{Proposition}
\emph{Proof.} 
First we recall that $\deg K =2g-2$ and $\rk|K|=g-1$. Since $\deg D=\deg K =2g-2$, we have that 
$$\deg (K-D)=\deg K-\deg D=0.$$
From Riemann-Roch, it follows
$$\rk|K-D| = \rk|D| - \deg D +g-1 \geq (g-1) -(2g-2)+g-1 = 0,$$
so $K \sim D$ and $|K|=|D|$.
\qed

\begin{Corollary} \label{iper}
On a tropical hyperelliptic curve $\G$ of genus $g \geq 2$, $$|K|=(g-1) \ g_2^1.$$
\end{Corollary} 
\emph{Proof.} 
Since $\G$ is hyperelliptic, there exists a (unique) $g_2^1$. 
If we define $|D|:=(g-1) \ g_2^1$ and use Proposition \ref{sistema}, we have finished, because $|D|$ has degree $2g-2$ and rank at least $g-1$.
\qed

\section{Clifford's Theorem}

Now we are in the position of proving the tropical version of Clifford's inequality:

\begin{Theorem}[Tropical Clifford] \label{clifford}
Let $D$ a special divisor on a tropical curve $\G$ of genus $g$. Then
$$\rk |D| \leq \frac{\deg D}{2}.$$
Furthermore, if $\G$ is hyperelliptic, equality occurs if and only if $D$ is a multiple of the (unique) $g_2^1$ on $\G$. 
\end{Theorem}
\emph{Proof.} 
First of all, as in the classical case, from the tropical Riemann-Roch Theorem for $D=K$, it follows that $\rk|K|=g-1$.

If $D$ is effective and special, then $K-D$ is also effective, so we can apply Lemma \ref{somma} and obtain
$$\rk|D|+\rk|K-D| \ \leq \ \rk|K| = g-1.$$
On the other hand, again by tropical Riemann-Roch, we have 
$$\rk|D|-\rk|K-D| \ = \ \deg D -g+1.$$
By adding these two expressions, we obtain the first part of the theorem:
$$2 \ \rk |D| \leq \deg D.$$

\bs

Assume now that $\G$ is hyperelliptic.

If $D$ is a multiple of the (unique) $g_2^1$, then obviously equality holds.

\ss

Conversely, we assume that equality holds and we 
consider the linear system 
$$|D| + (g-1-\rk|D|) \ g_2^1.$$
It has degree 
$$\deg D + (g-1-\rk|D|) \cdot 2 = 2g-2$$ 
and rank at least $g-1$ by Lemma \ref{somma}: indeed,
\begin{eqnarray*}
\rk \left(|D| + (g-1-\rk|D|) \ g_2^1 \right) &\geq& \rk |D| + \rk \left( (g-1-\rk|D|) \ g_2^1 \right) = \\
&=& \rk|D| + (g-1-\rk|D|) = g-1.
\end{eqnarray*}
Hence Proposition \ref{sistema} implies that it is the canonical system $|K|$, so by Corollary \ref{iper} we have
$$|D| + (g-1-\rk|D|) \ g_2^1 \ = \ |K| \ = \ (g-1) \ g_2^1.$$
Thus we conclude that $|D|=\rk|D| \ g_2^1$, which completes the proof. 
\qed

\bs

Every tropical curve of genus $g=2$ is hyperelliptic (see for instance \cite{BN2}, Example 5.1).
If $\G$ is a tropical curve of genus $g \leq 4$ and $D$ is an effective divisor on $\G$ such that
\begin{equation} \label{grado}
0 < \deg D < 2g-2
\end{equation}
and equality holds in Clifford's Theorem:
\begin{equation} \label{uguaglianza}
\rk|D|= \frac{\deg D}{2}
\end{equation}
then $\G$ is hyperelliptic. Indeed, if $g=3$, then from (\ref{grado}) and (\ref{uguaglianza}) it follows that $\deg D=2$ and $\rk|D|=1$.
If instead $g=4$ and $\deg D \neq 2$, then we have $\deg D = 4$ and $\rk|D|=2$, hence $\deg(K-D)=2$ and $\rk|K-D| = \rk|D|-\deg D + g -1 = 1$.

On the other hand, if $\G$ is a tropical curve of genus $g \geq 5$, then equality in Clifford's Theorem for a special divisor $D \neq 0,K$ on $\G$ does not imply \emph{a priori} that $\G$ is hyperelliptic.

\providecommand{\bysame}{\leavevmode\hbox to3em{\hrulefill}\thinspace}

\end{document}